Dmytro Taranovsky
December 1, 2015
Modified and Expanded: December 17, 2016

# Arithmetic with Limited Exponentiation

**Abstract:** We present and analyze a natural hierarchy of weak theories, develop analysis in them, and show that they are interpretable in bounded quantifier arithmetic $I\Delta_0$ (and hence in Robinson arithmetic Q). The strongest theories include computation corresponding to $k$-fold exponential (fixed $k$) time, Weak König's Lemma, and an arbitrary but fixed number of higher level function types with extensionality, recursive comprehension, and quantifier-free axiom of choice. We also explain why interpretability in $I\Delta_0$ is so rich, and how to get below it.

**Contents:**

**Outline:** In section 1, we explain that alternations of unbounded quantifiers is what makes interpretability in $I\Delta_0$ so rich, and that we can go below that by using logic with restricted alternations of unbounded quantifiers. Next (section 2), we discuss a relatively simple theory (using numbers, words, and predicates) and include relations with computational complexity classes. Then (section 3), we introduce the main theory, or rather a hierarchy of theories, using numbers, functions, and operators, and with extensions to higher types. Our choices aim to get the most in power, expressiveness, and convenience while staying within bounded quantifier definitions for the base theory and within interpretability in $I\Delta_0$ for all extensions. In section 4, we outline analysis and reverse mathematics in the theories, including a treatment of continuous functions without recursive comprehension and an equivalence for the intermediate value theorem. Finally (section 5), we prove that the theories are all interpretable in $I\Delta_0$ (and hence in Q).

# 1 Introduction and the Power of Quantifiers

## 1.1 The Hidden Strength of Unbounded Quantification

Part of foundations of mathematics is finding weak systems that capture much of mathematics, as well as understanding the relations between different systems and theorems. Exponential Function Arithmetic (EFA) captures most arithmetic, but can we go further by limiting exponentiation?

Formally, we can define weaker systems, such as $I\Delta_0$ (IDelta_0), which is arithmetic with bounded quantifier induction. However, we encounter surprising behavior.
**Definition:** $a_b^c$ is $a\hat{\ }a\hat{\ }\ldots\hat{\ }a\hat{\ }c$ ($b$ exponentials).
$I\Delta_0$ has length $n^{O(1)}$ proof that $2_n^2$ exists. More generally, given a model, a cut is a nonempty initial segment of the model that is closed under successor.
For every $k$, $I\Delta_0$ proves that there is a cut $I$ such that:
* $I$ is closed under '+' and '·', and
* for every $i \in I$, $2_k^i$ exists.
Thus, even in a very weak system, we can in a sense iterate exponentiation any fixed finite number of times. The proof is simple:
Set $I_0 = \mathbb{N}$, set $I_{n+1}(x) \Leftrightarrow (I_n$ is closed under $y \to y + 2^x)$. To get closure under '+' and '·', choose $\{n : \exists m \in I_{k+2}\, n < 2^{2^m}\}$.

The issue is not the induction — $I\Delta_0$ is interpretable Robinson Arithmetic Q, which has just a few of the properties of successor, addition, and multiplication. (For a good background, see "Interpretability in Robinson's Q" [Ferreira 2013].) Nor is it the totality of addition and multiplication; Q is interpretable in its weaker variant where we only require closure under successor [Švejdar 2008].

Instead, the strength comes from alternation of unbounded quantifiers. Meaningfulness of quantifiers over natural numbers implies a certain ability to do an unbounded search over natural numbers, and that $\mathbb{N}$ (the set of natural numbers) has certain attributes of a completed totality. The combinatorial strength of this principle is the partial ability to do exponentiation. Another way to view it is that unbounded quantifiers gives us access to higher types, and each higher type gives an exponential speed up. For example, set $F_0(f) = f * f$, $F_{n+1} = F_n * F_n$ ('*' is composition), so $F_n(f) = f^{2^{2^n}}$. This speed up characterizes the inherent strength of a basic ability to use higher types.

## 1.2 Restricted Quantifier Theories

Ordinary mathematics uses quantifiers in a limited way, and in particular uses only a very small number of unbounded quantifier alternations. By restricting quantifiers, we may restrict exponentiation while still capturing much mathematical practice. Specifically, given a notion of bounded quantifier or $\Sigma_0$ formulas, and a fixed small $k$, we can restrict the language to boolean combinations of $\Sigma_k$ formulas. There are different formalizations, but the following formalization may be best in terms of nonrestrictiveness and its closure properties:
- A $\Sigma_{i+1}$ formula is of the form $\exists x_1 \exists x_2 \ldots \exists x_n \varphi$ where φ is a boolean combination of

$\Sigma_i$ formulas.
- For the deduction system, use sequent calculus limited to boolean combinations of $\Sigma_k$ formulas.

Using cut elimination and the subformula property of cut-free proofs, and assuming that all axioms remain expressible, provability in the original theory agrees with the restricted theory for expressible formulas. The original theory may have an iterated exponential speed-up of proofs, but (at least if $k$ is not too small) that speed-up should only be likely for sentences that (perhaps indirectly) involve very large numbers.

**Notes:**
* One downside is that different reasonable proof systems might no longer be related up to polynomial size equivalence, but given a reasonable system $S$ for $k$ unbounded quantifiers and a reasonable system $T$ for $k+1$ unbounded quantifiers, $S$-proofs should be transformable with polynomial expansion into $T$ proofs. (Both claims are conjectural.)
* Intuitively, in the epistemic system corresponding to the proof system, the extent of natural numbers is fuzzy, and for $\Sigma_k$ $\varphi$, $\varphi(n)$ is meaningful for every concrete natural number $n$, but because of the fuzziness, the meaningfulness of $\forall n \varphi(n)$ is unclear.
* A different formalization is to limit validity to boolean combinations of $\Sigma_k$ sentences, and use a proof system based on sentences rather than open formulas. Note that a $k-1$ quantifier formula may correspond to a $k$ quantifier sentence.

Using cut elimination, the consistency of $I\Delta_0$ restricted this way is provable in EFA. Moreover, $I\Delta_0$ restricted to boolean combination of $\Pi_k^0$ formulas does not have length $n$ proof of existence of $2_{k+5}^n$.

Regarding expressiveness of restricted quantifier languages, for many statements k=1 suffices — most number theory theorems are $\Pi_1^0$ or have natural $\Pi_1^0$ strengthenings. However, it is annoying to always modify statements this way, and $k=2$ allows one to speak of (for example) an algorithm having linear time without specifying a proportionality constant. For second order theories (using predicates), basic axioms on real numbers appear to need at least $k=3$, and $k=4$ allows general $\Pi_2^1$ statements (and $k=5 - \Pi_3^1$). $\Pi_3^1$ statements are useful for exceptional theorems, $\Pi_2^1$ formulas, flexibility, and for allowing an extra quantifier for formulation in a weak base theory.

An alternative is to allow arbitrary formulas but limit the cut rule to restricted formulas, thus treating unrestricted formulas as having limited epistemic validity. This should get systems that — provably in EFA + iterated exponentiation — are equivalent to (for example) $I\Delta_0$, yet with consistency provable in EFA. However, absence of modus ponens is problematic, and allowing the cut rule for arbitrary sentences (but not open formulas) would allow length $n$ existence proof of about $2_{\log^*(n)}^n$. (One gets a cut-free proof that $I_k$ is closed under successor, and does cut elimination on $I_k(0), I_k(0) \Rightarrow I_k(1), I_k(1) \Rightarrow I_k(2), \ldots, I_k(n-1) \Rightarrow I_k(n)$.)

Let us define the *exponential index* of a theory as maximum $k$ such that existence of $2_k^n$ has a proof of length $n^{O(1)}$.

**Challenge:** Find natural theories with small exponential indices that capture much of mathematics.

If $S$ is a natural theory of exponential index $k$, then I conjecture that provably in a weak base theory, if $2_{k+2}^P$ exists, then $S$ has no inconsistency (i.e. a natural number coding an $S$-proof of $0 = 1$) less than $P$.

Besides exponential index, theories vary by the amount of induction, and one can intuitively associate a computational complexity class with a theory based on exponential index and the complexity corresponding to permitted induction. For theories with multiple types, there may be an ambiguity about which types to consider as numbers. By combining exponential index with the amount of induction we may ensure that we get the same complexity regardless of whether we use (say) unary or binary numbers. However, I do not have a formalization of this concept.

### 1.3 Importance of Interpretability in $I\Delta_0$

Setting the exact number of quantifiers may be arbitrary. Instead, we may step back and simply ask that theory has a finite exponential index if the number of unbounded quantifiers is fixed. However, a theory may show strength in ways other than exponential index (for example, through $\Pi_1^0$ strength), and the right formal measure is that the theory is interpretable in $I\Delta_0$. Moreover, we want there to be a formula in $I\Delta_0$ interpreting satisfaction for bounded quantifier formulas (whose length is in some definable cut) so that for restricted quantification we would not need to worry about such formulas.

It is surprising how much can be interpreted in $I\Delta_0$ — including Weak König's Lemma and much of mathematical analysis.

## 2 A Weak Theory and its Variations

### 2.1 Axioms and Discussion

Infinite objects can be interpreted using either predicates or functions. Using predicates leads to simpler axioms. However, in the absence of recursive comprehension, functions give much more flexibility. Functions are also more natural work with directly. In this section, we will give a simpler theory using predicates, and in the next section we will switch to using functions.

Even for finite programs, it is often useful to talk of input-time-output relation without noting implementation, and the natural language for this is to use infinite sets. Also, restricted exponentiation prevents equivalence of binary and unary numbers, so for maximum usability we want three types:
1. ordinary or unary numbers — as in physical sizes; closed under '$+$' and '$\cdot$'. One extension is to require closure under quasipolynomial functions, specifically $n \to \lfloor n^{\log n + 1} \rfloor$ ($n > 0$), and extend bounded quantifier comprehension accordingly.
2. words or binary numbers — as in data that a program can manipulate.
3. sets or predicates on words.

**Language:** Standard, except that we restrict to boolean combinations of $\Sigma_k$ formulas using $k = 3$ (see above). We can also use a different or infinite $k$.
**Axiomatization:**
* Basic axioms.
* Induction (PIND): Every nonempty set has a word of minimum length.
* Bounded quantifier comprehension:
$\forall P\, \exists Q\, \forall w\, (Q(w) \Leftrightarrow \varphi(P, w))$ where $\varphi$ is a bounded quantifier formula ($P$ may be a tuple of variables; $w$ is a word).
**Note:** The choice of "basic axioms" is unimportant and does not affect the theory as long as we can derive $S_2^1$ for binary numbers and that unary numbers are those whose exponent exists as binary. Also, a bounded quantifier is one of $\exists s < t$, $\exists v \subset t$, $\forall s < t$, $\forall v \subset t$, where $t$ is a unary term ($t$ does not use $v$ or $s$; $t$ may use '+' and '·'), $s$ is unary, and $v$ is binary (so $v \subset t \Leftrightarrow v < 2^t$).

Its exponential index is at least 5: Define $I_0 = \{x: 2_3^x$ exists as a binary number$\}$, $I_1 = \{x: I_0$ is closed under $y \to y + 2^x\}$, and $I_2 = \{x: I_1$ is closed under $y \to y + 2^x\}$. Membership in $I_2$ is $\Pi_3^0$ and hence meaningful for constants, allowing us to prove say $I_2(100)$. The exponential index rises to 9 if unary numbers are closed under quasipolynomial functions, and the system is extended to make sets with bounded word length count as numbers and behave accordingly, and we allow sentences with five blocks of unbounded quantifiers. (5 and 9 are lower bounds; I do not have a proof of optimality.)

One extension is to allow recursive comprehension. Even Weak König's Lemma can be added, but then we would lose optional interpretability of predicates as equivalence classes of programs. We will defer discussion to a later section. There is a hierarchy of systems based on how many times exponentiation can be effectively applied (and on other restrictions on comprehension and induction). A good project (partially accomplished here) would be to catalog some of these systems of arithmetic and analysis, show that their consistency is provable in EFA, and see how much mathematics can be done in them.

## 2.2 Relation with Complexity Classes

Bounded quantifier formulas correspond to the polynomial hierarchy PH, or relativized PH if the formula has a free set variable. There are natural comprehension schemes that (assuming basic comprehension) correspond to different complexity classes:
*Note:* Here, φ is bounded quantifier, and we can even pick a single polynomial time computable φ. PSPACE and EXPTIME are reducible to the formula shown; counting hierarchy can be reached level by level by using the previous level as a free variable. P is a free second order variable. The comprehension schema is that $Q$ is a set (or for counting hierarchy, a numeric function coded by a set).
* counting hierarchy: $Q(w) = |\{v < w : \varphi(P, v)\}|$
* PSPACE: $Q(w) = \varphi(P, \{v : Q(v) \wedge w - \text{len}(w) < v < w\})$
* EXPTIME: $Q(w) = \varphi(P, \{v < w : Q(v)\})$
* exponential hierarchy: Treat quantifiers over sets containing only words of restricted length as bounded.
* In the other direction (corresponding to $S_2^1$), one can restrict comprehension to

polynomial time P, but allow induction for NP formulas.

For the second order theory, the inclusions are proper (except possibly counting hierarchy v PSPACE) because there are oracle separations between the classes. However, for the first order theory, the only inclusion known to be proper is PH v. exponential hierarchy.

Capturing the counting hierarchy is useful to define acceptance probability of a polynomial time randomized (or even quantum) computation. Acceptance probability (after rescaling to get an integer) is complete for #P, which is conjectured to be outside PH; however approximations are in PH. For quantum computation, even approximations are not known to be in PH. Capturing #P also corresponds with the ability to do integration (under certain assumptions).

# 3 The Main Theory and its Extensions

## 3.1 Axioms and Basic Properties

The theory will have 4 types:
* unary numbers, as in physical sizes; this type is used for convenience and intuitive appeal; $n$ is unary iff $2^n$ exists as binary
* binary numbers, as in digital data; below, by numbers we will mean binary numbers
* functions: map numbers into numbers
* operators: map functions into functions.

**Notes:**
* The axioms will not mention unary numbers but they exist in a definitional extension by letting binary $2^n - 1$ code unary $n$ (or one can use a different coding as long as unary $n \to$ binary $2^n$ is bounded quantifier). One can also additional types of numbers, such as small numbers — say $n$ is small iff $2^n$ exists as unary (again with bounded-quantifier coding).
* Many $\Pi_2^0$ number theory theorems will need to weakened as: For every unary number there is a binary number such that the bounded-quantifier relation holds.
* One can add hyperoperators (mapping operators into operators) and in general any fixed finite number of higher operator types. Throughout the paper, we will indicate (sometimes briefly) how to extend the axioms and proofs to higher types.
* Because pairing will be codable, arity of functions and operators does not matter.
* By default, $m$ and $n$ will be numbers, $f$ and $g$ functions, and $F$ and $G$ operators.
* Equality is primitive only for numbers. $f = g$ means $\forall n f(n) = g(n)$, so $f = g$ is not a bounded quantifier formula.
* All functions and operators are total.
* The language is very expressive. For example, under the set-theoretical interpretation, the Continuum Hypothesis is equivalent to
$\exists F \forall f, g \exists n \, F(n, f) = g \lor F(n, g) = f$ ($n$ is a natural number and $f$ and $g$ are functions).

Define $\mu(x < n : \varphi(x))$ as the least number $x < n$ such that $\varphi(x)$ holds, and $n$ if there is no such number.

**Axioms:**

**1.** Basic arithmetical axioms. Include constants 0 and 1, and functions '+', '·', and $x \to x^{\text{len}(x)}$ (len is bit length; $\text{len}(x) = \lfloor \lg x \rfloor + 1$, except that $\text{len}(0) = 0$).
**2.** Induction for bounded quantifier formulas. Equivalently,
$\forall n\, f(n) = f(n+1) \Rightarrow \forall n\, f(n) = f(0)$.
**3.** Extensionality: $\forall n\, f(n) = g(n) \Rightarrow \forall m\, F(f)(m) = F(g)(m)$
**4.** Bounded quantifier comprehension:
   **a.** functions: $\exists f\, \forall n\, f(n) = \mu(m < T(n) : \varphi(m, n))$ ($T$ is a term; $T$ and $\varphi$ do not use $f$)
   **b.** operators: $\exists F\, \forall f\, \forall n\, F(f)(n) = \mu(m < T(n) : \varphi(m, n, f))$ ($T$ is a term; $T$ and $\varphi$ do not use $F$)
   where $\varphi$ is a bounded quantifier formula
   - each numeric quantifier is in the form $\forall z < S$ or $\exists z < S$ where $S$ is a term that does not use $z$.
   - each function quantifier is in the form 4a (i.e. one can use bounded minimization to define new functions).
   - each operator quantifier is in the form 4b.

**Notes:**
* Per convention, the axioms are open formulas that are assumed to be universally quantified.
* Bounded quantifier formulas appear very expressive: They can use functions and operators, and have quantifiers of the form $\exists n < f(m)$ (and $\exists n < F(f)(m)$ and so on), and even use $\mu$ to get new functions as input to operators.
* An alternative treatment would be to use a finite set of basic combinators (which would be operators for functions, and hyperoperators for operators) that capture bounded-quantifier comprehension and add their basic defining axioms. In that set up, bounded-quantifier definitions correspond to terms.
* If we used predicates instead of functions, and coded functions by predicates, then recursive comprehension would be necessary for basic things such as existence of $\{n : f(n) = 0\}$. By using functions directly, we avoid the need for unlimited recursive comprehension to get a meaningful theory of functions.

One shorthand is to use λ in place of bounded function and operator quantifiers. For example, set $\max = \lambda f.\, \lambda n.\, f(\text{argmax}(f)(n))$ where
$\text{argmax} = \lambda f.\, \lambda n.\, \mu(m \leq n : \forall m' \leq n\, f(m') \leq f(m))$, and we can show that $\max(f)(n) = \max_{m \leq n} f(m)$.

On the one hand, bounded quantifier comprehension suffices to give a basic general theory of functions and operators. On the other hand, the higher order theory with bounded quantifier comprehension is conservative over the first order theory with bounded-quantifier induction $I\Delta_0 + \Omega_1$. To see this, code functions and operators by sequences of bounded quantifier definitions, and note that evaluation allows one to expand them into numeric quantifiers and bounded quantifier formulas. However, this result is model-theoretic (i.e. the coding is done outside of the theory), and may involve a speed-up of proofs, especially for restricted quantification proofs. In particular, the definitions allow an exponential (the number of exponentials corresponding to the number of higher types) shortening of some bounded quantifier formulas.

Without recursive comprehension, the complexity of definitions of functions is characterized by two parameters:
(1) the complexity of evaluating the corresponding predicate
(2) the permitted growth rate.
High (2) leads to high (1). However, in a weak theory, (2) has to be very limited, so (1) is important to allow definition of complex functions. Without total exponentiation, even with recursive comprehension, one can define more complex predicates through bounded quantification over higher types (axiom 7 below).

*Polynomial Time Restriction:* To get a theory corresponding to polynomial time:
(1) Add unary numbers (with exponentiation function to get binary numbers).
(2) Treat ordinary bounded quantifiers over binary numbers as unbounded.
(3) In place of μ, use $\sum_{i=0}^{\text{len}(T)} 2^i \varphi(i)$ ($i$ is unary and $\varphi$ is bounded quantifier; 0 is False, and 1 is True) (and treat the corresponding definitional quantifier for a binary number as bounded).
(4) Add $\Sigma_1^b$-PIND: $\varphi(0) \wedge \forall n(\varphi(\lfloor n/2 \rfloor) \Rightarrow \varphi(n)) \Rightarrow \forall n \varphi(n)$, where $\varphi$ is $\exists m \leq n\, p(m,n)$ and $p$ is a boolean function symbol.
(5) Bounded Time Comprehension: Treat the following quantifier as bounded: $\exists w \forall i\, (w(i) \Leftrightarrow i < j \wedge q(w \upharpoonright i))$ ($w$ is a word (a binary number), $i$ and $j$ are unary ($j$ is a free variable), and $q$ is a boolean function).
*Note:* It may be more natural to treat unary numbers as numbers here and binary numbers as words.
*Extension:* To compensate for the lack of max operator, we can add $\forall n \exists n' \forall m \leq n\, f(m) < n'$ or the stronger $\forall f \exists g \forall n \forall m \leq n\, f(m) < g(n)$.
*Strength:* I did not verify whether the above theory is conservative over $S_2^1$.

## 3.2 Extensions

**5.** Bounded collection: $\forall m \exists n\, f(m,n) = 1 \Rightarrow \forall m \exists n \forall m' < m \exists n' < n\, f(m', n') = 1$
**Note:** Bounded collection (even for polynomial time relation $f$) appears unprovable in PRA, even if extended with all true $\Pi_2^0$ statements, but it is very useful for dealing with arbitrary recursive functions (without recursive comprehension, recursive 'functions' need not be functions in the theory).

**6.** Add (or assert existence of) SUM operator: $\text{SUM}(f)(n) = \sum_{i=0}^{n} f(n)$.
**7.** More generally, treat quantification over $P(n)$ (using functions) or $P(P(n))$ (using operators) as bounded, or otherwise allow a high computational complexity class in defining functions and operators. (Here, $P$ stands for power set, and $n$ is a number. A boolean function $f$ codes an element of $P(n)$ iff $\forall m \geq n\, f(m) = 0$, and a boolean $F$ codes an element of $P(P(n))$ iff $F(f)$ depends on $f$ only through $\{m < n : f(m) = 0\}$.)
**8.** Allow quasipolynomial time functions: $x \to x^{\text{len}(x)^{\text{len}(\text{len}(x))}}$.

For a fixed $k$ (but not simultaneously for all $k$ in $\mathbb{N}$), we can even allow
$\Omega_k$: For all $x$, x^len(x)^...^len$^k$(x) exists, where len$^k$ is $k$-fold iteration of len.
Function $\Omega_k$: x→x^len(x)^...^len$^k$(x) exists as a function.

*Extension to higher types:* The language and the theory can be extended with hyperoperators (and any fixed finite number of types) defined analogously to operators, and satisfying extensionality and bounded quantifier comprehension, and axiom 7 can be extended to treat quantifiers over $P(P(P(n)))$ (and any fixed number of power sets) as bounded. One formalization of bounded quantification over $P(\ldots P(n)..)$ is to treat the quantifier $\exists G(G = C_n(G))$ as bounded where $G$ is a function (or a higher type), and $C_n$ 'contracts' its argument to $n$, as $C$ does in the proof of axiom 12.

*Integers of higher types:* In the absence of total exponentiation, and using axioms 1-7, we can still iterate exponentiation of ordinary numbers a fixed number of times using numbers of higher types (and with exponentiation being a function/operator /etc). For example, we can treat $f: \mathbb{N} \to \{0,1\}$ with $f(0) = \sup(x+1 : f(x+1) = 1)$ as $\sum_{i=1}^{f(0)} 2^{i-1} f(i)$ to get numbers in $2^{<\mathbb{N}} = \{n : \text{len}(n) \in \mathbb{N}\}$; let $\text{code}(n)$ refer to this coding of numbers by functions. Similarly, we can use operators to get numbers in $2^{<2^{<\mathbb{N}}} = \{n : \text{len}(\text{len}(n)) \in \mathbb{N}\}$, and so on. The exact coding does not matter as long as basic operations are bounded quantifier. Note that while there is no bounded quantifier formula testing whether a function codes a number (i.e whether $f \in \text{code}$), there is an operator $F$ such that $F(f) \in \text{code}$ and also $f \in \text{code} \Rightarrow F(f) = f$, which suffices. For example, the quantifier $\exists G : 2^{<\mathbb{N}} \to 2^{<\mathbb{N}}$ is equivalent to the quantifier $\exists G' \exists G = \lambda f. F(G'(F(f)))$ and hence uses only one unbounded quantifier. We may also want to add bounded collection as axioms for numbers of higher types if we do not include stronger axioms that imply it.

**Further Extensions:**
**9.** Recursive Comprehension: $\forall m \exists n \, p(m,n) \Rightarrow \exists f \forall m \, p(m, f(m))$ ($p$ is a boolean function).
**10.** Weak König's Lemma (WKL): Every infinite binary tree has an infinite path. Here, a tree can be coded by a boolean $f$ with $\forall n > 0 \, f(n) \Rightarrow f(\lfloor (n-1)/2 \rfloor)$ ($\lfloor (n-1)/2 \rfloor$ is the parent of $n$), and a path is a nonbranching tree.

**Notes:**
* Axioms 9 and 10 are extremely useful for mathematical analysis, but 9 breaks the relation between functions and bounded time computability, and 10 breaks the relation between functions and computability.
* WKL implies recursive comprehension (and $\Sigma_1^0$ separation) for predicates on unary numbers.
* There is an extension of WKL that we can add:
**10a.** WKL holds for numbers of the type $2^{<\mathbb{N}}$ (coded by functions, and with the tree and path coded by operators) (and its extensions to higher types).
In the absence of exponentiation, 10a does not appear to follow from WKL.
* If exponentiation is not total, then in the presence of 1-7, WKL is implied by (code is as defined in "Integers of higher types" (above)):
**10b.** $\forall f : \mathbb{N} \to \{0,1\} \, \exists F \, \forall n \in \mathbb{N} \, F(\text{code}(n)) = f(n)$
The implication holds because (using $F$) we can evaluate $f$ past $\mathbb{N}$, and process the binary tree coded by $f$. Axiom 7 ensures that numbers in $2^{<\mathbb{N}}$ are well-behaved. It is unclear if the statement is provable from the axioms, but it (including extensions to $2^{<\mathbb{N}}$ and so on) is interpretable in $I\Delta_0$ (using the second proof in "Proof of Interpretability in $I\Delta_0$" section).

**Further extensions for higher types:**
**11.** Recursive Comprehension for Operators: $\forall f \exists n \, P(f, n) \Rightarrow \exists F \forall f \, P(f, F(f)(0))$ ($P$ is a boolean operator) (and similar extensions for operators of higher types).
**12.** Uniformization (or quantifier-free axiom of choice) (extension of 11):
$\forall f \exists g \, P(f, g) \Rightarrow \exists F \forall f \, P(f, F(f))$ ($P$ is a boolean operator) (and similar extensions for operators of higher types).

**Notes:**
* There may be different reasonable notions of recursive comprehension for operators and higher types, and it would be interesting to work out different theories, their strength, and relations with extensionality and WKL.
* Axiom 11 appears the best match for recursive comprehension (when used in conjunction with bounded quantifier comprehension) while retaining interpretability in $I\Delta_0$.
* Axiom 12 looks very powerful, and under the set semantics, it is unprovable in ZF. However, the power depends on certain prerequisites that are not available in the base theory, and in the model and $F$ we will construct, $F(f)$ will be (in the appropriate sense) finite. Also, in the model, every operator $G$ that is a predicate on boolean functions is finite. By contrast, existence of an operator testing for a boolean function being zero would lead to Peano Arithmetic (PA). Also, one weakening of 12 provable in ZF is $\forall f \exists! g \exists h \, P(f, g, h) \Rightarrow \exists F \forall f \exists h \, P(f, F(f), h)$, but (unless $g$ is a number) the premise (or just $\exists P \exists! g \, P(g)$) implies PA (and $ACA_0$ for predicates).

In place of a hierarchy of theories, we can also speak of their union (thus using infinite number of types, bounded quantification over $P(P \dots P(N) \dots)$, and closure under $\Omega_k$). The resulting theory will have every finite fragment interpretable in $I\Delta_0$. It appears that the combined system proves the same $\Pi_1^0$ statements as EFA, but with an iterated exponential slow-down of certain proofs. (The amalgamation of numbers of higher types forms a model of EFA, but we cannot get to $n$ exponentials using less than $n$ symbols.)

### 3.3 Formula Expressiveness

A $\Pi_k$ sentence is a sentence of the type $QX_1 QX_2 \dots QX_k \, \varphi$, where each $Q$ is a quantifier ($\forall$ or $\exists$), the first $Q$ being $\forall$ (or $\exists$ for $\Sigma_k$), each $X_i$ is a tuple of variables, and $\varphi$ is a bounded quantifier formula.
Under the set-theoretical interpretation, $\Pi_{k+1}$ formulas correspond to $\Pi_k^2$ formulas ($k > 0$) (and with hyperoperators, $\Pi_k^3$ formulas, and so on). A $\Sigma_1^2$ formula can be coded as $\exists F \forall g \, \varphi$ where $\varphi$ is bounded quantifier. (It can also be coded as $\exists F \forall G \, \varphi$, where $\varphi$ is an equation and $F$ is binary.)
Conversely (assuming comprehension, which is not in our axioms), a $\Sigma_1$ formula can be witnessed by an operator which is zero on all but a finite set of functions. $\Sigma_1$ statements using functions and operators are essentially $\Sigma_1^0$, and all true $\Sigma_1$ statements are provable from the axioms 1-4 (even if $\varphi$ may quantify over $P(n)$ and $P(P(n))$). The only complication is extensionality, which is addressed by giving a witness that the functions on which the operators gave different results are in fact different. I conjecture that this (provability of true $\Sigma_1$ statements in the analog of axioms 1-4) applies to extensions to higher types and that a true $\Sigma_1$ statement has a finite (for example, as in Finite Example Property in the proof of axiom 12) example.

# 4 Real Numbers and Analysis

## 4.1 Introduction and Real Numbers

Reverse mathematics up to $RCA_0$ illuminated much of the structure of mathematical analysis. Going further, and doing analysis without full recursive comprehension will illuminate which constructions use bounded versus unbounded time. Studying weak theories also shows the complexity of certain constructions, such as the equivalence between (a) being able to do integration, (b) having SUM operator (axiom 6), and (c) having oracles for relativized PP problems.

Because of absence of exponentiation, there are two natural notions of real numbers, both of which satisfy real field axioms. For every unary number $t$, ordinary real numbers are defined with $2^{-t}$ error, and approximate real numbers are defined with $1/t$ error. They can can defined as follows.
* Ordinary real numbers: Let $(n, f)$ code $n + \sum_{i=1}^{\infty} f(2^i)/2^i$ where $n$ is an integer and each $f(2^i)$ is in $\{0, 1, 2\}$.
* Approximate real numbers: Same as above, but $n$ is a unary integer and use $f(2^{2^i})$ in place of $f(2^i)$.
Operations between ordinary and approximate real numbers yield approximate real numbers, except when the result is out of bounds. Approximate real numbers are useful especially for weaker theories, such as with polynomial time computation. In the other direction, given enough comprehension, one can also define real numbers corresponding to integers of higher types.

A function $\mathbb{R} \to \mathbb{R}$ may depend only the value of the real number (that is $r = s \Rightarrow f(r) = f(s)$ even though $f(r)$ and $f(s)$ may be different as codes), as opposed to an intensional function that may depend on the code.

## 4.2 Continuous Functions

Without recursive comprehension and WKL, classically equivalent notions become nonequivalent. For continuity, the differences are in (1) to what extent must a modulus of continuity exist as a function (and included in the code, the later being relevant for operators on functions and quantifier counting), and (2) whether the continuity holds for (in a sense) generic real numbers (or points) or just the ones in the model. One approach is to just require (and include in the code) a modulus of pointwise continuity as a function. Instead, however, we will identify a general category of spaces and impose the needed conditions to get a good theory with just axioms 1-4.

*Type of space:* Complete metric spaces with every bounded subset totally bounded. Examples include $\mathbb{R}^n$ (unary $n$) and the Hilbert cube (each successive dimension halving in size).
*Coding of the Space:* Such a space is coded by using a countable dense set (indexed by numbers), a metric $d$ (on that dense set), and a function witnessing total boundness: $h(r, \varepsilon) \to m$ such that every point within distance $r$ from the first point is within distance $\varepsilon$ of one of the first $m$ points; for coding, $r$ and $1/\varepsilon$ can be restricted to (and coded as) numbers.
*Coding of points:* Analogous to real numbers. For example, use a sequence $x$ (of indices in the dense set) such that $d(x_n, x_{n+1}) < 2^{-n}$.
*Coding of a function $f$:* (1) Specify (as a function of the index) the value of $f$ for every element of the dense set. (Alternatively, use any reasonable standard coding of continuous functions using functions rather than predicates.)
(2) Give a function witnessing that $f$ is uniformly continuous on every bounded set: $g(r, \delta) \to \varepsilon$ such that $d(y, 0) < r \wedge d(y, z) < \varepsilon \Rightarrow d(f(y), f(z)) < \delta$ and $\delta > 0 \Rightarrow \varepsilon > 0$, with $r, 1/\delta, 1/\varepsilon$ restricted to and coded by positive integers.
*Generalization to other spaces (not formalized):* Given a function $h$ with for all $x$, $h(x)$ being a compact subset of the domain $f$, require there to be $g$ with $g(x)$ returning a modulus of uniform continuity for $f$ on $h(x)$. If in our domain there is a canonical choice of $h$, include $g$ in the code for $f$.

With this definition, continuous functions are coded by numeric functions rather than operators, function application is bounded quantifier (with invalid arguments potentially returning invalid values), and the choice of coding is irrelevant (in various senses that can be made precise). Quantification over functions uses a single unbounded quantifier: There is an operator that changes every code $f$ into a code $f'$ such that $f'$ is valid, and if $f$ is valid, $f$ and $f'$ code the same function. The requirement of uniform continuity on bounded subsets (which automatically holds under recursive comprehension and WKL) allows us to develop the intended basic theory in a weak system. Also, for typical continuous functions on $[0, 1]^n$, the theory suffices to prove continuity, and for typical continuous $f : \mathbb{R}^n \to \mathbb{R}^n$, the theory suffices unless $f$ grows or changes at a superpolynomial rate.

Note that by having $g$ as a function, we do not need recursive comprehension for basic theory: For example, the supremum (equivalently a code for the supremum) of a continuous function $f$ on $[0, 1]^n$ can be computed by a bounded quantifier formula (using $g$). Also, using axioms 1-4,6, we can show that integration on $\mathbb{R}^n$ exists as an operator, and is linear and satisfies most of its basic properties; axiom 6 is not needed for approximate reals (if unary numbers are closed under $x^n$).

### 4.3 Reverse Mathematics

We can also do reverse mathematics.

**Theorem (axioms 1-5):** Intermediate value theorem is equivalent to recursive comprehension for predicates on unary numbers.
**Proof:** ($\Leftarrow$) If the function has no rational zeros, use the constructive version (below) and wait until the interval becomes short enough before returning it. If the interval does not become shorter than $\varepsilon$, then the function is 0 on an interval of length $\varepsilon$ and hence has a rational zero.
($\Rightarrow$) Given a recursive (allowing parameters) predicate $g$, we can find a bounded quantifier (allowing parameters) function $f$ on reals with one zero, and that zero will code $g$. The construction is as follows. Fix an interval. Run an evaluation of $g(0)$, and during that time, output that $f$ is close to 0 on that interval and is negative to the left of the interval and positive to the right of the interval. Then, if $g(0) = 0$, narrow the interval to the first third, and if $g(0) = 1$ — the last third, and continue with $g(1)$ and so on.
**Notes:**
* In the base system, the complexity of a function is decoupled from its growth rate, so finding complex predicates does not prove existence of fast growing functions. Also, complex predicates on unary numbers need not imply complex nonsparse predicates on binary numbers. However, given bounded collection (for binary numbers), recursive comprehension for functions on unary numbers implies recursive comprehension for functions on binary numbers.
* Without bounded collection, the equivalence is for this form of comprehension: $\forall m \exists n \forall m' \leq m \exists n' \leq n \, p(m', n') \Rightarrow \exists g \, \forall i \, (g(i) \Leftrightarrow \min(n : p(2^i, n))$ is odd) ($g$ is a predicate on unary numbers, and $p$ is a predicate).

**Problem:** Find other natural equivalents to recursive comprehension.

A constructive weakening of intermediate value theorem is provable in the base system: If $f$ is continuous, and $f(0) < 0$ and $f(1) > 0$, then there is $g$ such that $g(n)$ returns $(a_n, b_n)$ with $\forall x (a_n < x < b_n) \, |f(x)| < 1/n$ and $a_j < a_k < b_k < a_j$ for $j < k$.

We can show that every continuous function on $[0, 1]$ has a supremum. However, having a maximum value, along with many other statements in analysis, is equivalent to Weak König's Lemma. Some of these equivalences have been worked out in a theory called BTFA and its extensions (see [Ferreira 1994] and [Fernandes]). *Subsystems of Second Order Arithmetic* [Simpson 2009] is a good reference for reverse mathematics, including equivalences of Weak König's Lemma assuming $RCA_0$. (Note that in these references, continuity on [0,1] is not defined so as to give uniform continuity, and hence in their treatment $RCA_0$ does not prove that continuous functions on $[0, 1]$ are bounded and have other basic properties.) Many of these equivalences apply here assuming axioms 1-6,9.

# 5 Proof of Interpretability in $I\Delta_0$

## 5.1 Setup

Despite its expressiveness and apparent power, the theory (axioms 1-12 and extensions) is interpretable in $I\Delta_0$. Furthermore, under our interpretation in $I\Delta_0$, there is a definable cut $K$ and a formula $\varphi$ coding the satisfaction relation for bounded quantifier formulas whose length is in $K$.

In $I\Delta_0$ (and even in extended Robinson arithmetic $Q^+$), we have a rich collection of definable cuts:
* Given a definable cut $I$, there is a definable cut $J$ that such that for $\forall n \in J \, 2^n \in I$
* Unless $I$ satisfies $\Sigma_2$ induction (and $\Sigma_k$ induction), there is a definable cut $J < I$.

**Proposition:** A model of $I\Delta_0$ has no definable cuts iff it satisfies Peano Arithmetic.
**Proof:** A definable cut contradicts induction. Conversely, given an induction schema $\forall n \, \varphi(n)$, consider $\{n : \forall m \leq n \, \varphi(m)\}$. If the set is a cut, it is a definable cut, otherwise choose minimum $n$ that is not in the set, and use the induction instance coded by $\varphi(n)$ to get a cut.

If there is too much induction ($\Sigma_1$ induction suffices), then instead of getting a proper cut of $\mathbb{N}$, we will have (in the construction below) $J = \mathbb{N}$ and use a coding for nonstandard numbers to get $I$. Alternatively, known results on interpretability of WKL$_0$ might be extendable to our axioms. Technically, the interpretation will behave one way if sufficient induction holds and another way if it fails.

In its general form, the theory can be parameterized by
$k_1$ — numbers satisfy $\Omega_{k_1}$
$k_2$ — the number of higher types (i.e. $k_2 = 2$ if we use just functions and operators).
$k_3$ — number of exponentials sufficient for computation of bounded quantification (with axiom 7 and its analogues $k_3 = k_2 + 1$, and without it $k_3 = 1$).
$k_4$ — (not used in the proof) for restricted quantification theories, maximum $k$ such that $\Pi_k$ formulas are allowed.

We present (or sketch) two proofs/constructions, one in which operators essentially correspond to bounded quantifier definitions, and another in which the set of operators is in a sense maximized. The first proof only applies to axioms 1-10 (not 10a), and in that proof, I did not verify extensionality for hyperoperators and higher types.

## 5.2 First Proof

The construction will use three cuts: $K \leq J < I$.
$I$ — satisfies $\Omega_1$
$J$ — satisfies $\Omega_{k_1}$, and with $\exists n \notin J \, 2^n_{k_3} \in I$.
$K$ — closed under '+' and '·', and such that for every total function $f$ on $J$ (as below, with a number coding the graph of $f$ up to some point beyond $J$): $J$ is also closed under $f^{2^n_{k_2}}$ for all $n$ in $K$.

**Interpretation**
\* Numbers are elements of $J$.
\* Functions are arbitrary functions that are total on $J$ with a number coding the graph of the function up to some point beyond $J$.
\* Primitive operators (and similarly for higher types): Every number $x$ with len(len($x$))∈J is included as an operator constant (and similarly for higher types). The exact coding is unimportant as long as the codes can be processed with bounded quantifier formulas.
\* Operators (and hyperoperators, etc.) and bounded quantifier formulas are coded by a sequence of bounded-quantifier definitions whose total length (number of bits) is in $K$, and that may refer to numbers and functions and primitive operators (and if applicable, hyperoperators, etc).

**Formula Evaluation**
*Goal:* Given a bounded quantifier formula and values for its free variables (number, functions, primitive operators, function and operator definitions), evaluate whether it is true.
*Procedure:* The formula is evaluated recursively, keeping track of free variables.
\* $A \vee B$ is evaluated by evaluating $A$ and evaluating $B$, and similarly with $A \wedge B$.

* $\exists m < T$ is evaluated by evaluating $T$, and then evaluating the formula for every $m < T$,
and similarly with $\forall m < T$ and $\mu(m < T : \varphi)$, and also for quantification over $P(n)$, etc.
* For definitional quantifiers, note the definition and proceed to evaluate the formula.
* It remains to show how to evaluate a term. Each top level term is in one of three forms:
  - $n$ where $n$ is a number
  - $f(T)$ where $f$ is a function symbol. Evaluate $T$. Then evaluate $f(T)$ either directly or recursively using the definition of $f$.
  - $F(T_f)(T)$ where $F$ is an operator symbol and $T_f$ is a function term. Evaluate $T$, add function $T_f$, and evaluate $F(T_f)(T)$ using the definition of $F$. Note that for the primitive operators above, $F(f)$ can be evaluated by looking up its value on (under one choice of coding) $(f(0), f(1), \ldots, f(f(0)))$.
  - Higher types are handled analogously.

**Notes:**
* Bounded collection holds because its failure would allow computation of the least element not in $J$.
* Recursive comprehension holds because given $f$ (coded as $f : n \to n$ for some $n \notin J$), we can evaluate a recursive-in-$f$ function on $J$ up to some element beyond $J$, and halt when the computation reaches $n$, which will only happen outside of $J$.
* WKL can be shown similarly. A coded predicate $P$ on $J$ is defined up to some point beyond $J$, and we can choose a path through the binary tree corresponding to $P$.
* Closure of $K$ under multiplication ensures that given a formula and substitutions for some of its free variables, the result is valid.
* Closure of $I$ and $J$ under the relevant functions ensures that bounded quantifier comprehension using formulas with length in K still leads to valid numbers and functions. It also ensures that formula evaluation (above) is bounded-quantifier (in $I$), assuming that we are given a sufficiently large element of $I$ (independent of the formula) and may return an arbitrary result if a parameter is invalid or the formula length is not in $K$.
* A length $n$ bounded quantifier formula can apply a function up to about $2^n_{k_2}$ number of times. Conversely, for an operator $F$ with bounded-quantifier definition length $n$, we can find a function $g$ such that $F(f)(m)$ uses $f$ only for arguments less than $h^{2^n_{k_2}}(m)$ where $h(x) = \max_{y \leq x}(f(y) + g(y))$. Similar results also hold for higher types, and extensionality follows. In $I\Delta_0$, while we cannot do induction on extensionality directly, we can use induction for these concrete properties.
* Note that the same function can have multiple codes; it is not clear if this can be avoided (and it would be interesting to get non-interpretability results under unique coding).
* We can vary the construction to get (in a sense) minimal models for various subsystems. For examples, for axioms 1-4, we can simply code functions and operators by bounded quantifier definitions.
* The interpretation (and also the interpretation in the second proof) satisfies two desirable properties: It is definable parameter free, and by being cut-based, $I\Delta_0$ proves that for every $\Sigma^0_1$ formula whose length is in K the interpretation satisfies only true instances.

## 5.3 Second Proof

Let $I$ and $J$ (and numbers and functions) be as above. We ensure extensionality while keeping function (including higher level function) application bounded quantifier (in $I$) as follows.

**Coding:** Objects of type $k$ (0 is numbers, 1 is functions, and so on) will be coded by valid level $k$ codes. Given $n$ ($n \notin J$ and $2^n_{k_2} \in I$), a level 0 code is a number $< n$; a level $k+1$ code is a function (coded by a number): $\text{code}_k \to \text{code}_k$ where $\text{code}_k$ is the set of all level $k$ codes. A level 0 code is valid iff it is in $J$; equivalence is identity. A level $k+1$ code is valid iff it maps equivalent valid level $k$ codes into equivalent valid level $k$ codes. Two valid level $k+1$ codes $f$ and $g$ are equivalent iff $f(x)$ is equivalent to $g(x)$ for every valid level $k$ code $x$.

**Note:** This construction appears ideal for getting the most while staying within interpretability in $I\Delta_0$. A good project may be to investigate whether there are additional natural true axioms that this construction can be made to satisfy and are interpretable in $I\Delta_0$.

**Formula Evaluation:** When encountering a definitional quantifier, compute the table corresponding to the function (or operator, etc.) defined in the quantifier. (We get a point outside of $J$ using one of the higher-order parameters, or if there are none, we can pick (or quantify over) any point that is not too large.)

**Cut $K$:** Let $K \subset J$ be the cut (or a subcut closed under '·') such that bounded quantifier formulas whose length is in $K$ are extensional, that is their values do not depend on how the parameters are coded, and such that terms whose length is in $K$ evaluate to valid numbers (or functions, etc.). We have to show that this is a cut. Extensionality is clearly preserved by logical connectives and bounded numerical quantifiers.

**Definitional Quantifiers:** Extensionality is preserved by the definitional quantifiers assuming that the definition actually (and extensionally) defines a function (or operator, etc.). A level $k$ function $f_k$ (level 2 being operators) can be characterized as an extensional function $f: f_{k-1}, \ldots, f_2, f_1, m \to n$, and we need to show that a bounded quantifier definition (as in the axioms) guarantees existence and extensionality of $f$. Existence follows from formula evaluation and closure of $J$ under terms whose length is in $K$. Extensionality follows from extensionality of the underlying formula for (essentially) $f(f_{k-1}, \ldots, f_2, f_1, m) = n$.

**Proof of Axiom 12:**

*Extensional Extensions and Contractions:* If in the codes (for functions, operators, etc.), we switch coding from using $N$ (and $N \to N$, etc.) to using $N'$ (and $N' \to N'$, etc.), we want to extend (if $N' > N$) or contract (if $N' < N$) the functions (including higher level functions) while preserving extensionality on $J < \min(N', N)$. Here is one recipe, where $E$ is extension and $C$ is contraction:
$E(n) = n$, $C(n) = \min(\min(N', N) - 1, n)$, $E(F)(f) = E(F(C(f)))$, $C(F)(f) = C(F(E(f)))$.
Note that as desired, $C(E(F))(f) = C(E(F)(E(f))) = C(E(F(C(E(f))))) = F(f)$.
One can also show (by recursion on the type level) that extensionality is preserved (but note that contractions may add extensionality), and that for valid $F$ and $f$, $C(F)(C(f))$ is equivalent to $C(F(f))$, and also that extension from an $N$ in $J$ to an $N'$ outside of $J$ guarantees extensionality.

*Finite Example Property:* If $\exists F\, P(F)$ ($F$ and $P$ are valid), then there is a finite '$F$' (finite means being based on some $N \in J$) whose extension satisfies $P$.

*Proof:* Using extensional contractions, a suitable '$F$' exists for every $N \notin J$, and since $J$ is a cut, it must exist for an $N \in J$.

*Axiom 12 Proof:* Assuming $\forall f \exists g\, P(f, g)$ ($P$ is a boolean operator/hyperoperator /etc.), set $F(f)$ to be $E(g)$ for the least (which implies using the least number $N$) $g$ such that $P(f, E(g))$.

*End of proof*

Proofs of other axioms are routine. Proof of WKL is as in the first proof.

## 5.4 Further Remarks

**Continuity** *(not part of the proof)*: It appears that in the model (in the second proof), higher level functions are precisely the hereditarily continuous extensional functionals (coded by numbers in $I$). A code for $F$ (with $F(f) \in \mathbb{N}$) is a set of partial_input-output pairs (each pair is a number of type $n - 1$ where the type of $F$ is $n$) such that for every code $c$ for $f$, every sufficiently large initial segment of $c$ will give $F(f)$. Note that codes are non-unique, and extensionality is not vacuous. The property of being a code for $F$ is definable (and code application is partial recursive and total on valid codes), and we could have developed the theory using codes at the price of requiring the continuity, imposing properties such as recursive comprehension on the base theory, and having a cumbersome exposition, including lack of a simple analog to our bounded quantifier formulas.

Using cut elimination (see for example [Beckmann 2011] and [Gerhardy 2005]), the interpretability implies that for every finite $k_1, k_2, k_3, k_4$ if we restrict first order logic to boolean combinations of $\Pi_{k_4}$ formulas, there is $k$ such that every proof of existence of $2^n_k$ has length $> n$.

## REFERENCES

Beckmann, Arnold; Buss, Samuel R. Corrected upper bounds for free-cut elimination. [J] Theor. Comput. Sci. 412, No. 39, 5433-5445 (2011). ISSN 0304-3975.

Fernandes, António M.; Ferreira, Fernando. Basic applications of weak König's lemma in feasible analysis.

Ferreira, Fernando. A Feasible Theory for Analysis. J. Symbolic Logic 59 (1994), no. 3, 1001--1011.

Ferreira, Fernando; Ferreira; Gilda. Interpretability in Robinson's Q. Bulletin of Symbolic Logic, 19, pp 289-317 (2013). doi:10.1017/S1079898600010660.

Gerhardy, Philipp. The Role of Quantifier Alternations in Cut Elimination. Notre Dame J. Formal Logic, Volume 46, Number 2 (2005), 165-171.

Simpson, Stephen G. Subsystems of Second Order Arithmetic (Perspectives in Logic) 2nd Edition (2009).

Švejdar, Vítězslav. An Interpretation of Robinson Arithmetic in its Grzegorczyk's Weaker Variant. Fundamenta Informaticae. Volume 81 Issue 1-3, May 2008, 347-354.